\documentstyle[amscd,amssymb,verbatim,epsf,mathrsfs,graphicx,hyperref]{amsart}

\begin{document}
\theoremstyle{plain}
\newtheorem{Thm}{Theorem}
\newtheorem{Cor}{Corollary}
\newtheorem{Ex}{Example}
\newtheorem{Con}{Conjecture}
\newtheorem{Main}{Main Theorem}
\newtheorem{Lem}{Lemma}
\newtheorem{Prop}{Proposition}

\theoremstyle{definition}
\newtheorem{Def}{Definition}
\newtheorem{Note}{Note}

\theoremstyle{remark}
\newtheorem{notation}{Notation}
\renewcommand{\thenotation}{}

\errorcontextlines=0
\numberwithin{equation}{section}
\renewcommand{\rm}{\normalshape}%

\title[Global Joachimsthal] {A Global Version of a Classical Result of Joachimsthal}
\author{Brendan Guilfoyle}
\address{Brendan Guilfoyle\\
          School of Science, Technology, Engineering and Mathematics \\
          Institute of Technology, Tralee \\
          Clash \\
          Tralee  \\
          Co. Kerry \\
          Ireland.}
\email{brendan.guilfoyle@@ittralee.ie}
\author{Wilhelm Klingenberg}
\address{Wilhelm Klingenberg\\
 Department of Mathematical Sciences\\
 University of Durham\\
 Durham DH1 3LE\\
 United Kingdom.}
\email{wilhelm.klingenberg@@durham.ac.uk }
\date{18th February 2016}


\begin{abstract}
A classical result attributed to Joachimsthal in 1846 states that if two surfaces intersect with constant angle along a line of curvature of one surface, 
then the curve of intersection is also a line of curvature of the other surface. In this note we prove a global analogue of this result, as follows. 

Suppose that two closed convex surfaces intersect with constant angle 
along a curve that is not umbilic in either surface. We prove that the principal foliations 
of the two surfaces along the curve are either both orientable, or both non-orientable.

We prove this by characterizing the constant angle intersection of two surfaces in Euclidean 3-space as the intersection of a surface 
and a hypersurface in the space of oriented lines. The surface is Lagrangian, while the hypersurface is null, with respect to the canonical 
neutral Kaehler structure. We establish a relationship between the principal directions of the two surfaces along the intersection curve in Euclidean 
space, which yields the result.

This method of proof is motivated by topology and, in particular, the slice problem for curves in the boundary of a 4-manifold.
\end{abstract}

\maketitle

\section{Introduction}

A classical result of differential geometry states that if two surfaces, $S_1$ and $S_2$, in ${\mathbb R}^3$ intersect with constant angle along a line 
of curvature of $S_1$, then the curve of intersection is also a line of curvature of $S_2$.  In 1846 Joachimsthal proved a special case of this 
theorem \cite{Joach}, while the general case was presented by Bonnet in 1853 \cite{Bonnet}. 

In this note we prove a global analogue of this result. That is, consider two closed convex surfaces $S_1$ and $S_2$ intersecting transversely along a 
smooth
simple closed curve $C$. Decompose $S_1=D_1^+\cup_CD_1^-$, $S_2=D_2^+\cup_CD_2^-$, where $D_j^\pm$ are closed discs, with $D_1^-$ defined to be 
the disc with outward pointing normal given by the normal to $S_2$ projected onto $S_1$, and similarly for $D_2^-$. We prove

\vspace{0.1in}
\begin{quote}
\noindent{\bf Main Theorem.}
{\it If $S_1$ and $S_2$ intersect with constant angle along a curve that is not umbilic in either $S_1$ or $S_2$, then the principal foliations 
of the two surfaces along the curve are either both orientable, or both non-orientable.}
\end{quote}
\vspace{0.1in}

The Main Theorem is proven by characterizing the constant angle intersection of surfaces in ${\mathbb R}^3$ as the intersection of a surface $\Sigma$ 
and a hypersurface ${\cal H}_\epsilon$ in the space ${\mathbb L}({\mathbb R}^3)$ of oriented lines of ${\mathbb R}^3$. Here $\epsilon=\tan(\alpha/2)$, 
where $\alpha$ is the angle of intersection. In this characterization, the surface $\Sigma$ is a Lagrangian  section and ${\cal H}_\epsilon $ is a 
null hypersurface with respect to the neutral K\"ahler structure on the space of oriented lines. From these geometric properties we 
relate the principal directions of the surfaces along the intersection (equation (\ref{e:final})) and the result follows.

Our methods are motivated by topology: the ultimate goal is to investigate whether a knot in the boundary of a 4-manifold is the boundary of a 
properly embedded disc - the slice problem. As a new geometric model, we propose subsets of ${\mathbb L}({\mathbb R}^3)$ 
and invariants derived from the neutral K\"ahler structure, which we claim can be utilized to detect obstructions to sliceness. 

For surfaces in ${\mathbb R}^3$, umbilic points are generically isolated and have a half-integer index associated with them given by the winding number
of the principal foliation about the umbilic point. Let $D$ be a disc with non-umbilic boundary and only isolated umbilic points in the interior.
The total umbilic index $i(D)$ is defined to be the sum of the umbilic indices on $D$ and when the boundary of $D$ is non-umbilic, it is equal
to the winding number of the principal foliation along $\partial D$. 

From our perspective, umbilic points on $S$ correspond to complex points on $\Sigma$ and the total
complex index $I(\Sigma)=2i(D)$ of the corresponding disc in ${\mathbb L}({\mathbb R}^3)$ with totally real boundary is a smooth topological 
invariant \cite{Lai}. Thus our main Theorem can be restated as:
\[
I(\Sigma_1^+)=I(\Sigma_2^+) {\mbox{ mod 2}}
\] 
Alternatively, our result provides a situation where one can bound this invariant by geometric data at the boundary hypersurface. Further results on 
this topic and its relation to the slice problem will be reported in a forthcoming paper.

\vspace{0.2in}

\section{Proof of the Main Theorem}

Let $S_1$ and $S_2$ be smooth convex surfaces which intersect along a curve $C$ at a constant angle $0<\alpha<\pi$.  

Split the convex surfaces so that $S_1=D_1^+\cup_CD_1^-$, $S_2=D_2^+\cup_CD_2^-$, where $D_j^\pm$ are closed discs, with $D_1^-$ defined to be 
the disc with outward pointing normal given by the normal to $S_2$ projected onto $S_1$, and similarly for $D_2^-$, as in the diagram below.
\vspace{0.1in}

\setlength{\epsfxsize}{3.0in}
\begin{center}
{\mbox{\epsfbox{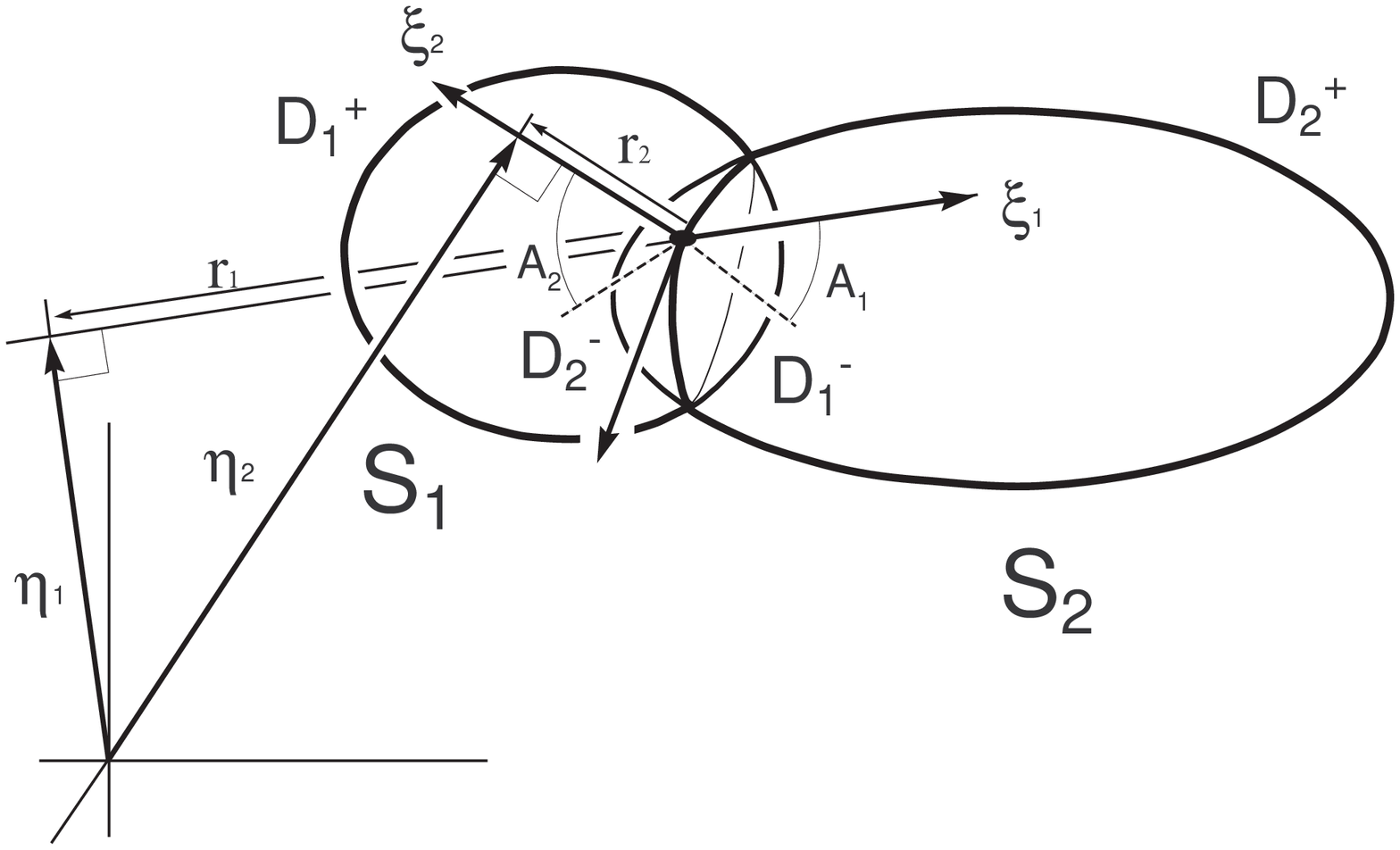}}}
\end{center}

\vspace{0.1in}

Now transfer this geometric data to the space  ${\mathbb L}({\mathbb R}^3)$ of oriented lines of Euclidean ${\mathbb R}^3$. ${\mathbb L}({\mathbb R}^3)$
admits a canonical neutral K\"ahler structure $({\mathbb J},{\mathbb G},\Omega)$ that is invariant under the induced Euclidean action \cite{gak4}, and, 
since  ${\mathbb L}({\mathbb R}^3)$ is diffeomorphic to $TS^2$, write $T{\mathbb S}^2$ for ${\mathbb L}({\mathbb R}^3)$ endowed with 
this structure.

The surfaces $S_1$ and $S_2$ in ${\mathbb R}^3$ will be treated differently: define ${\cal H}_\epsilon \subset T{\mathbb S}^2$ to be the set 
oriented lines that intersect $S_1$, forming an angle $\alpha$ with the outward normal, and define $\Sigma\subset T{\mathbb S}^2$ be the 
set of oriented lines normal to $S_2$.
Clearly these are 3- and 2-dimensional submanifolds of $T{\mathbb S}^2$, the latter being diffeomorphic to a 2-sphere while the former is 
diffeomorphic to a circle bundle over the 2-sphere . 

The submanifolds ${\cal H}_\epsilon $ and $\Sigma$ intersect transversely along the curve ${\cal C}\subset T{\mathbb S}^2$ which consists of the 
oriented normal lines to $S_2$ along the curve $C=S_1\cap S_2\subset {\mathbb R}^3$. The constant angle condition in 
${\mathbb R}^3$ has been translated to an intersection condition for certain submanifolds in $T{\mathbb S}^2$. 
We now describe the geometric properties of these submanifolds. 

Starting with $\Sigma$, the set of oriented normals to the convex surface
$S_2$. This is a Lagrangian section of the bundle $T{\mathbb S}^2\rightarrow {\mathbb S}^2$ and is locally described by the gradient of the support 
function $r_2:{\mathbb S}^2\rightarrow {\mathbb R}$. 

We now utilize local holomorphic coordinates $(\xi,\eta)$ on $T{\mathbb S}^2$ to describe this  - see \cite{gak4} and references therein for
further details. The complex coordinate $\xi$ is the standard holomorphic coordinate on 
${\mathbb S}^2-\{{\mbox {south pole}}\}$ giving the direction of the oriented line, while the complex coordinate $\eta$ gives the perpendicular 
distance of the line to the origin. Together they yield 
\[
(\xi,\eta)\leftrightarrow \eta\frac{\partial}{\partial\xi}+\bar{\eta}\frac{\partial}{\partial\bar{\xi}}\in T_\xi{\mathbb S}^2.
\]
The surface $\Sigma$ is the section $\xi\mapsto(\xi,\eta=\eta_2(\xi,\bar{\xi}))$
where 
\begin{equation}\label{e:pot1}
\eta_2={\textstyle{\frac{1}{2}}}(1+\xi\bar{\xi})^2\frac{\partial r_2}{\partial \xi},
\end{equation}
and this is Lagrangian as
\[
\frac{\partial}{\partial \xi}\left(\frac{\eta_2}{(1+\xi\bar{\xi})^2}\right)
=\frac{\partial}{\partial \bar{\xi}}\left(\frac{\bar{\eta}_2}{(1+\xi\bar{\xi})^2}\right).
\]
Define
\[
\sigma_2=-\frac{\partial \bar{\eta}_2}{\partial \xi}
\qquad\qquad
\psi_2=r_2+(1+\xi\bar{\xi})^2\frac{\partial}{\partial \xi}\left(\frac{\eta_2}{(1+\xi\bar{\xi})^2}\right)
\qquad\qquad
\kappa_2=\psi_2^2-|\sigma_2|^2.
\]
If $R_1$ and $R_2$ are the radii of curvature of $S_2$, then 
\[
\sigma_2\bar{\sigma}_2={\textstyle{\frac{1}{4}}}(R_1-R_2)^2
\qquad\qquad
\psi_2={\textstyle{\frac{1}{2}}}(R_1+R_2)
\qquad\qquad
\kappa_2=R_1R_2.
\]
Thus umbilic points on $S_2$ are precisely the complex points on $\Sigma$, while the half of the argument of $\sigma_2$ determines the principal 
directions of the surface $S_2$. 

Let $(\xi,\eta=\eta_1(\xi,\bar{\xi}))$, $r_1$, $\sigma_1$, 
$\psi_1$ and $\kappa_1$ be the analogous quantities for the surface $S_1$ and its associated Lagrangian section.

We turn now to ${\cal H}_\epsilon $, the oriented lines making a fixed angle $\alpha$ with $S_1$. To write this in local coordinates, note the 
following: 

\vspace{0.1in}
\begin{Lem}\label{l:1}
Given $\xi_1\in {\mathbb S}^2$, the points on ${\mathbb S}^2$ forming an angle $0<\alpha<\pi$ with $\xi_1$ is a circle parameterized by
\[
\xi_2=\frac{\xi_1+\epsilon e^{iA_1}}{1-\bar{\xi}_1\epsilon e^{iA_1}},
\]
for $A_1\in[0,2\pi)$, where $\epsilon=\tan(\alpha/2)$.
\end{Lem}
\begin{pf}
If $\xi_1\leftrightarrow\vec{X}_1$ and $\xi_2\leftrightarrow\vec{X}_2$, then
\[
\vec{X}_1\cdot\vec{X}_2=\frac{2\xi_1\bar{\xi}_2+2\xi_2\bar{\xi}_1+(1-\xi_1\bar{\xi}_1)(1-\xi_2\bar{\xi}_2)}{(1+\xi_1\bar{\xi}_1)(1+\xi_2\bar{\xi}_2)}
   =\frac{1-\epsilon^2}{1+\epsilon^2}=\cos\alpha,
\]
as claimed.
\end{pf}
\vspace{0.1in}

Note that the above relationship can be inverted to
\[
\xi_1=\frac{\xi_2+\epsilon e^{iA_2}}{1-\bar{\xi}_2\epsilon e^{iA_2}},
\]
where the angles $A_1$ and $A_2$ are related by
\begin{equation}\label{e:a1a2}
e^{iA_2}=-\frac{e^{iA_1}-\epsilon\xi_1}{e^{-iA_1}-\epsilon\bar{\xi}_1}e^{-iA_1}.
\end{equation}

Now suppose we have an oriented line $(\xi_1,\eta_1)$ and a point $(z,t)$ on this line. Thus, for some $r_1\in{\mathbb R}$ 
\begin{equation}\label{e:coordu}
z=\frac{2(\eta_1-\bar{\xi}_1^2\bar{\eta}_1)}{(1+\xi_1\bar{\xi}_1)^2}+\frac{2\xi_1}{1+\xi_1\bar{\xi}_1}r_1
\qquad\qquad
t=-\frac{2(\xi_1\bar{\eta}_1+\bar{\xi}_1\eta_1)}{(1+\xi_1\bar{\xi}_1)^2}+\frac{1-\xi_1\bar{\xi}_1}{1+\xi_1\bar{\xi}_1}r_1.
\end{equation}
Suppose further that $(\xi,\eta)$ is another oriented line that passes through $(z,t)$ and forms an angle $\alpha$ with $(\xi_1,\eta_1)$. 
The former fact means that:
\[
\eta={\textstyle{\frac{1}{2}}}(z-2t\xi_2-\bar{z}\xi_2^2).
\]
Substituting the expressions for $z$ and $t$, and for $\xi_2$ from Lemma \ref{l:1}, we conclude that
the set of oriented lines that intersect $(\xi_1,\eta_1)$ at a parameter $r=r_1$, forming an angle $\alpha$, is a circle in $T{\mathbb S}^2$
parameterized by $A_1\in[0,2\pi)$ as follows
\[
\xi=\frac{\xi_1+\epsilon e^{iA_1}}{1-\bar{\xi}_1\epsilon e^{iA_1}}
\qquad\qquad
\eta=\frac{\eta_1-\epsilon^2e^{2iA_1}\bar{\eta}_1-\epsilon(1+\xi_1\bar{\xi}_1)e^{iA_1}r_1}{(1-\bar{\xi}_1\epsilon e^{iA_1})^2},
\]
where $\epsilon=\tan(\alpha/2)$. 

As $(\xi_1,\eta_1)$ vary over the surface set of normal lines to $S_1$, the above equations trace out the hypersurface 
${\cal H}_\epsilon \subset T{\mathbb S}^2$, which is in fact null with respect to the neutral K\"ahler metric.

By assumption $\Sigma$ and ${\cal H}_\epsilon$ intersect in the curve ${\cal C}$, consisting of the oriented 
normal lines to $S_2$ along the curve of intersection $C=S_1\cap S_2$. Parameterize the curve so that $C$ is given by equations 
(\ref{e:coordu}) where everything now depends upon a parameter $u$.

Note that
\begin{equation}\label{e:deta1}
\frac{d \eta_1}{du}=\frac{\partial\eta_1}{\partial\xi_1}\frac{d \xi_1}{du}+\frac{\partial\eta_1}{\partial\bar{\xi}_1}\frac{d\bar{\xi}_1}{du}
=\left(\psi_1-r_1+\frac{2\bar{\xi}_1\eta_1}{1+\xi_1\bar{\xi}_1}\right)\frac{d \xi_1}{du}-\bar{\sigma}_1\frac{d\bar{\xi}_1}{du},
\end{equation}
and
\begin{equation}\label{e:dr1}
\frac{d r_1}{du}=\frac{\partial r_1}{\partial\xi_1}\frac{d \xi_1}{du}+\frac{\partial r_1}{\partial\bar{\xi}_1}\frac{d\bar{\xi}_1}{du}
=\frac{2\bar{\eta}_1}{(1+\xi_1\bar{\xi}_1)^2}\frac{d \xi_1}{du}+\frac{2\eta_1}{(1+\xi_1\bar{\xi}_1)^2}\frac{d\bar{\xi}_1}{du},
\end{equation}
where we have used equation (\ref{e:pot1}) and the definitions of $\sigma_1$ and $\psi_1$. Identical expressions hold for the derivatives of 
$\eta_2$ and $r_2$ with all subscripts changed from 1 to 2.

This yields the following:

\vspace{0.1in}
\begin{Lem}
To fix unit parameterization of the curve $C$ set
\begin{equation}\label{e:unitpara}
\frac{d\xi_1}{du}=\frac{(1+\xi_1\bar{\xi}_1)}{2\kappa_1}\left[\psi_1e^{i(A_1+\beta)}+\bar{\sigma}_1e^{-i(A_1+\beta)}\right],
\end{equation}
for $\beta(u)\in[0,2\pi)$. The curve ${\cal C}$ of oriented lines are tangent to $C$ iff $\beta=0$ or $\pi$, while, the oriented lines are normal to 
$C$ iff $\beta=\pi/2$ or $3\pi/2$.
\end{Lem}
\vspace{0.1in}

In summary, the intersection curve ${\cal C}$ in $T{\mathbb S}^2$ is the solution to the two complex equations
\begin{equation}\label{e:int}
\xi_2=\frac{\xi_1+\epsilon e^{iA_1}}{1-\bar{\xi}_1\epsilon e^{iA_1}}
\qquad\qquad
\eta_2=\frac{\eta_1-\epsilon^2 e^{2iA_1}\bar{\eta}_1-(1+\xi_1\bar{\xi}_1)\epsilon e^{iA_1}r_1}{(1-\bar{\xi}_1\epsilon e^{iA_1})^2},
\end{equation}
where all quantities depend on $u$. In the second equation we use the given graph function $\eta_1(u)=\eta_1(\xi_1(u),\bar{\xi}_1(u))$, and 
similarly for $\eta_2$. Thus we have five real unknowns encoded in $(\xi_1,\xi_2,A_1)$ and four real equations, yielding a one dimensional solution set
parameterized by $u$. 

Now differentiate the first of these equations with respect to $u$ and use equation (\ref{e:unitpara}) with $\beta=pi/2$ to get
\begin{align}\label{e:xi2dot}
\frac{d \xi_2}{du}&=\frac{i(1+\xi_1\bar{\xi}_1)e^{iA_1}}{2(1-\bar{\xi}_1\epsilon e^{iA_1})^2}\left[
2\epsilon\frac{dA_1}{du}+\epsilon(\epsilon e^{2iA_1}+e^{iA_1}\xi_1)\frac{\sigma_1}{\kappa_1}\right.\nonumber\\
&\qquad+\left.(e^{-iA_1}\epsilon\bar{\xi}_1-e^{-2iA_1})\frac{\bar{\sigma}_1}{\kappa_1}
-(\epsilon(e^{iA_1}\bar{\xi}_1+e^{-iA_1}\xi_1)+\epsilon^2-1)\frac{\psi_1}{\kappa_1}\right].
\end{align}

Differentiate the second equation of (\ref{e:int}) with respect to $u$ and use equations  (\ref{e:dr1}), (\ref{e:unitpara}) 
and (\ref{e:xi2dot}), along with (\ref{e:deta1}) and the analogous equation for $\eta_2$. The result, after some computation and rearrangement, is
\begin{align}
&(1-e^{iA_1}\epsilon\bar{\xi}_1)^2\bar{\sigma}_2\Bigg[2\epsilon\kappa_1\frac{d A_1}{du}
   +(-e^{2iA_1}+e^{iA_1}\epsilon \xi_1)\sigma_1\nonumber\\
&\qquad\qquad\qquad+\epsilon(e^{-iA_1}\bar{\xi}_1+\epsilon e^{-2iA_1})\bar{\sigma}_1
   -(\epsilon(e^{iA_1}\bar{\xi}_1+e^{-iA_1}\xi_1)+\epsilon^2-1)\psi_1\Bigg]\nonumber\\
&+(e^{iA_1}-\epsilon\xi_1)^2\psi_2\Bigg[2\epsilon\kappa_1\frac{d A_1}{du}
   +\epsilon(\epsilon e^{2iA_1}+e^{iA_1}\xi_1)\sigma_1\nonumber\\
&\qquad\qquad\qquad+(e^{-iA_1}\epsilon\bar{\xi}_1- e^{-2iA_1})\bar{\sigma}_1
   -(\epsilon(e^{iA_1}\bar{\xi}_1+e^{-iA_1}\xi_1)+\epsilon^2-1)\psi_1\Bigg]\nonumber\\
&\qquad\qquad\qquad\qquad-(1+\epsilon^2)\kappa_1(e^{iA_1}-\epsilon\xi_1)^2=0.\nonumber
\end{align}

Taking the complex conjugate of this equation and eliminating the derivative of $A_1$ we get
\[
\frac{\sigma_1e^{2A_1i}-\bar{\sigma}_1e^{-2A_1i}}{\kappa_1}
   -\frac{\sigma_2}{\kappa_2}\frac{(e^{iA_1}-\epsilon\xi_1)^2}{(1-e^{iA_1}\epsilon\bar{\xi}_1)^2}
+\frac{\bar{\sigma}_2}{\kappa_2}\frac{(1-e^{iA_1}\epsilon\bar{\xi}_1)^2}{(e^{iA_1}-\epsilon\xi_1)^2}=0.
\]
Finally, recalling equation (\ref{e:a1a2}), we obtain
\[
\frac{\sigma_1e^{2A_1i}-\bar{\sigma}_1e^{-2A_1i}}{\kappa_1}
   -\frac{\sigma_2e^{2A_2i}-\bar{\sigma}_2e^{-2A_2i}}{\kappa_2}=0.
\]

To transfer this to back to Euclidean 3-space, denote the principal curvatures of $S_1$ by $(\lambda_1,\mu_1)$ and those of $S_2$ by $(\lambda_2,\mu_2)$. 
Let $\phi_1$ be the angle between the principal direction of $S_1$ and the intersection curve and $\phi_2$ the corresponding angle on $S_2$. 
The last equation can now be succinctly written
\begin{equation}\label{e:final}
(\lambda_1-\mu_1)\sin\phi_1-(\lambda_2-\mu_2)\sin\phi_2=0.
\end{equation}
Thus we obtain Joachimsthal's Theorem:  $C$ is a line of curvature on $S_1$ iff it is a line of curvature on $S_2$.
A further result of Joachimsthal's also follows: if $S_1$ is a sphere (so that $\sigma_1=0$) and the surface $S_2$ intersects it at constant angle,
then  $\sigma_2e^{2iA_2}-\bar{\sigma}_2e^{-2iA_2}=0$ and so $C$ is a line of curvature on $S_2$. Note that this equation states that the geodesic 
torsion of the intersection curve on the two surfaces are equal \cite{doC}.

Moreover, if $C$ contains no umbilic points on $S_1$ or $S_2$, the winding numbers of the maps $\phi_1,\phi_2:S^1\rightarrow S^1$ count the total 
umbilic index inside the curve on each surface, where the inside is determined by the outward pointing normals - see Theorem 4 of \cite{gak4}.

For simplicity, write $f_1=\lambda_1-\mu_1$ and $f_2=\lambda_2-\mu_2$ so that the derivative of equation (\ref{e:final}) with respect to $u$ is
\begin{equation}\label{e:las2}
f'_1\sin\phi_1+f_1\cos\phi_1\;\phi'_1=f'_2\sin\phi_2+f_2\cos\phi_2\;\phi'_2
\end{equation}
If $\phi_1\neq\phi_2$ everywhere along the curve then the winding numbers of the principal foliations are equal and the Theorem holds. Otherwise,
choose a point $p\in C=S_1\cap S_2$ such that $\phi_1(p)=\phi_2(p)$. By a deformation of the surfaces preserving constant angle we can also 
assure that $\phi'_1(p)\neq0$ (an open condition). 

By a choice of parameterization of the intersection curve $C$ we can set $p$ to be $u=0$ and $u=2\pi$. Furthermore, by a rotation of the complex 
coordinates $\xi$ we can set $\phi_1(0)=\phi_2(0)=0$. Then $\phi_1(2\pi)=n\pi$ and $\phi_2(2\pi)=m\pi$ for some $n,m\in{\mathbb Z}$. Each
principal foliation is orientable or non-orientable along $C$ according to whether the winding numbers $n$ or $m$ is even or odd, respectively.

Now evaluating equation (\ref{e:las2}) at $u=0$ and $u=2\pi$ we find
\[
f_1(p)\phi'_1(p)=f_2(p)\phi'_2(p) \qquad\qquad (-1)^nf_1(p)\phi'_1(p)=(-1)^mf_2(p)\phi'_2(p). 
\]
Since none of $f_1, f_2, \phi'_1$ or $\phi'_2$ vanish at $p$, both of these equations hold only if $n=m+2k$ for some $k\in{\mathbb Z}$. Thus the 
foliations are both orientable or both non-orientable, as claimed. 

\qed

\end{document}